\documentclass[runningheads]{llncs}

\usepackage{graphicx}
\usepackage{epstopdf}
\usepackage[caption=false]{subfig}

\usepackage{amsmath, amssymb}
\usepackage{mathtools}
\usepackage{graphicx}
\mathtoolsset{showonlyrefs}
\usepackage{xcolor}
\usepackage{accents}

\def\pd#1{{\color{blue}#1}} 
 
\def\dv#1{{\color{magenta}#1}}

\def\pd#1{{\color{black}#1}} 
\def\dv#1{{\color{black}#1}}

\def\pddelete#1{{\color{blue}#1}} 
\def\dvdelete#1{{\color{magenta}#1}}
\renewcommand{\pddelete}[1]{\unskip}
\renewcommand{\dvdelete}[1]{\unskip}

\usepackage{hyperref}
\hypersetup{
    colorlinks=true,
    citecolor=blue,
    linkcolor=red, 
    filecolor=magenta,      
    urlcolor=cyan,
}

\usepackage{wrapfig}

\usepackage{xcolor}
\usepackage{subfig}
\usepackage{hyperref}
\usepackage{textcomp}
\usepackage{nicefrac}

\usepackage{algorithm}
\usepackage{algpseudocode}
 \usepackage{algorithmicx}
\algnewcommand{\LeftComment}[1]{\Statex \(\triangleright\) #1}

\usepackage{xcolor}
\usepackage{bbm}

\usepackage{graphicx}
\begin{document}

\title{Stochastic optimization for dynamic pricing\thanks{The research is supported by the Ministry of Science and Higher Education of the Russian Federation (Goszadaniye) №075-00337-20-03, project No. 0714-2020-0005.}}

\author{Dmitry Pasechnyuk\inst{1}\orcidID{0000-0002-1208-1659} \and
Pavel Dvurechensky\inst{2}\orcidID{0000-0003-1201-2343} \and 
Sergey Omelchenko\inst{1}\orcidID{0000-0002-7266-3917} \and 
Alexander Gasnikov\inst{1}\orcidID{0000-0002-7386-039X}}

\authorrunning{D. Pasechnyuk et al.}

\institute{Moscow Institute of Physics and Technology, Russia
\and Weierstrass Institute for Applied Analysis and Stochastics, Germany}

\maketitle
\begin{abstract}
We consider the problem of supply and demand balancing that is stated as a minimization problem for the total expected revenue function describing the behavior of both consumers and suppliers. In the considered market model we assume that consumers follow the discrete choice demand model, while suppliers are equipped with some quantity adjustment costs. The resulting optimization problem is smooth and convex making it amenable for application of efficient optimization algorithms with the aim of automatically setting prices for online marketplaces. We propose to use stochastic gradient methods to solve the above problem. We interpret the stochastic oracle as a response to the behavior of a random market participant, consumer or supplier. This allows us to interpret the considered algorithms and describe a suitable behavior of consumers and suppliers that leads to fast convergence to the equilibrium in a close to the real marketplace environment.

\keywords{Automatic pricing \and Expected revenue function \and Stochastic convex optimization \and Supply and demand balancing.}
\end{abstract}

\section{Introduction}

With the development of platforms for online trading and services provision, the problem of dynamic pricing becomes more and more urgent. The environment of online marketplaces and financial services raises the question of establishing the most relevant prices for the items presented, and in such cases, relevance is understood as the possibility of market clearing, i.e. balancing supply and demand.

There are various approaches to finding equilibrium prices. In this paper, we consider an approach based on the characterization of the market state by a potential function of total expected revenue, similar to the function of total excessive revenue from \cite{nesterov2017distributed}. In \cite{muller2021dynamic}, the authors, using this approach, analyze the case in which consumers follow a discrete choice model with imperfect behavior introduced by random noise in the assessment of utility, and suppliers seek to maximize their profit, taking into account the quantity adjustment costs. This paper substantiates the convexity of the used potential as a function of prices. Hence, it follows that there are prices that minimize the potential, and it turns out that this minimum satisfies the condition of market clearing.

Thus, if there is an intermediary responsible for the formation of prices, it only needs to build a sequence of price values leading to a minimum. At the same time, it seeks to use the smallest possible amount of observations of market participants. Note that the introduced formulation is very convenient for applying the results of optimization theory to find the optimal prices. Indeed, the function under consideration is to be minimized and it simultaneously includes the characteristics of both consumers and suppliers, although in the classical game-theoretic approach the opposition of the interests of these parties leads through the Nash equilibrium scheme to minimax problems. Nevertheless, it turns out to be possible to propose exactly the potential function that describes the whole system and reduces its dynamic to only tendency to an extreme point, similar to physical systems. On the other hand, the considered function has the form of a sum including together the terms that characterize the behavior of both parties and each agent. This means that we can consider both consumers and suppliers as the same market participants and uniformly take into account their interests, without distinguishing between their types. Together, the described advantages lead to the ability to use the developments of convex optimization to find methods for constructing a sequence of prices, and analyze them in terms of the oracle complexity, which in this setting means the number of observations of the agents' behavior.

In this paper, continuing the ideas of \cite{muller2021dynamic}, the problem of finding equilibrium prices is posed in a stochastic setting, and the number of observations of the intermediary for single market participants is determined as a measure of the effectiveness of the method for solving this problem. Section~\ref{S:problem} describes the problem statement and introduces the stochastic oracle for considered function to optimize. Section~\ref{algorithms} discusses various algorithms for stochastic convex optimization, proposes theoretical estimates of the efficiency of the algorithms, and describes the practical advantages of each of them. Section~\ref{infty_d} generalizes the stochastic setting for the case of an infinite number of consumers, which makes it possible to take into account the previously unobserved consumers, and describes the algorithms for this setting. Section~\ref{zero_gamma} considers a special case of the problem with zero quantity of the adjustment costs, in which the potential loses its smoothness property, and proposes several approaches to get around this difficulty. Section~\ref{experiments} shows the results of modeling of the behavior of the proposed algorithms for a synthetic problem and demonstrates the improvement in the efficiency of the methods, which can be achieved by considering the stochastic formulation of the problem, in comparison with the dynamics from \cite{muller2021dynamic}.

\section{Problem statement}
\label{S:problem}

Let us imagine the marketplace environment. There are a number of suppliers (for example, shops) and a number of consumers (buyers). And every supplier offers some products (goods) in the assortment. Products are grouped by their type, every group contains some alternatives. In turn, every consumer can choose to buy one of the alternatives, guided by its subjective utility and price (with some element of randomness). For reasons of increasing profit, suppliers can change the assortment, taking into account the costs of these changes themselves. We use the corresponding mathematical models for the consumer and the supplier, and characterize them by the values of expected surplus and maximal revenue, respectively. Summing these terms up we obtain the function of total expected revenue. It describes the current imbalance in the market system. Its minimum corresponds to the prices at which the market is cleared. This means that we can formulate an optimization problem for finding equilibrium prices.

Let us consider such an optimization problem (in a simplified form) for the case of $n$ product alternatives divided into $m$ disjoint groups $G_i \subset \{1, ..., n\}$, $S$ suppliers with convex costs functions $c_s: \mathbb{R}^n_+ \rightarrow \mathbb{R}_+$, closed and convex \pddelete{boundedness is not required, right? }\dvdelete{Right. It can be bounded for economic reasons, but we do not use this property. In \cite{muller2021dynamic} it is not bounded too } sets of capacity constraints $Y_s \subset \mathbb{R}^n_+$, typical supplies $\hat{y}_s \in Y_s$ and quantity adjustment costs equal to $- \Gamma_s \cdot \|y - \hat{y}_s\|_2^2$ for some parameter $\Gamma_s > 0$, $D$ consumers with matrix $A = \{a_{i d} > 0\}_{1, 1}^{n, D}$ of alternatives subjective utility for each of them, and prices vector $p$:
\begin{equation}\label{finite_problem}
    \min_{p \in \mathbb{R}^n_+} \left\{ f(p) := \sum_{s=1}^S \pi_s(p) + \sum_{d=1}^D E_d(p) \right\},
\end{equation}
where \dvdelete{Clarified above. }\pddelete{What is $\Gamma_s$?}
\begin{equation}\label{pi}
    \pi_s(p) = \max_{y \in Y_s} \left\{ \langle y, p \rangle - c_s(y) - \Gamma_s \cdot \|y - \hat{y}_s\|_2^2 \right\}
\end{equation}
is the maximal revenue of supplier with taking into account costs $c_s(\cdot)$ and quantity adjustment costs parametrized by $\Gamma_s$ for given prices, and
\begin{equation}\label{E}
    E_d(p) = \ln  \left( \sum_{j=1}^m \left( \sum_{i \in G_j} e^{(a_{i d}-p_i)/\mu_j} \right)^{\mu_j}\right)
\end{equation}
is the expected surplus $\displaystyle E_d(p) = \mathbb{E}_\epsilon\left[\max_i \left\{a_{i d} - p_i+ \epsilon_i\right\}\right]$ for the discrete choice demand model with noise and for corresponding nested logit distribution \cite{mcfadden1978modeling}, with some correlation parameters $0 < \mu_j \leq 1$. \dvdelete{Solved: }\pddelete{I suggest to give more explanations on why the problem is stated in this way. I Understand (but still would be nice to mention this explicitly) that $\pi_s(p)$ is motivated by maximization of the revenue for given prices $p$. But, the meaning of the regularization term with $\Gamma_s$ is not clear. Also it is not clear, why $E_d(p)$ is constructed in this way and why the whole objective is to minimize in $p$.} We can also provide an explicit expression for the gradient of the introduced objective function:
\begin{equation}\label{grad}
    \nabla f(p) = \sum_{s=1}^S y_s(p) - \sum_{d=1}^D x_d(p),
\end{equation}
where $y_s(p)$ is the optimal solution of the optimization problem \eqref{pi} (by Demyanov–Danskin theorem \cite{danskin2012theory}), and
\[
[x_d(p)]_i = \frac{\displaystyle e^{(a_{i d} - p_i) / \mu_j} \left( \sum_{k \in G_j} e^{(a_{k d} - p_k) / \mu_j} \right)^{\mu_j - 1}}{\displaystyle  \sum_{h=1}^m \left(\sum_{k \in G_h} e^{(a_{k d} - p_k) / \mu_h}\right)^{\mu_h}},\quad i \in \{1, ..., n\}
\]
are the probabilities to choose certain alternative by the consumer\footnote{Since $\displaystyle E_d(p) = \mathbb{E}_\epsilon\left[\max_i \left\{a_{i d} - p_i+ \epsilon_i\right\}\right] = \sum_i \mathbb{P}\left[a_{i d} - p_i+ \epsilon_i = \max_i \left\{a_{i d} - p_i+ \epsilon_i\right\}\right] \cdot \mathbb{E}_\epsilon\left[a_{i d} - p_i+ \epsilon_i\right]$, and therefore $\displaystyle \frac{\partial E_d(p)}{\partial p_i} = -\mathbb{P}\left[a_{i d} - p_i+ \epsilon_i = \max_i \left\{a_{i d} - p_i+ \epsilon_i\right\}\right]$}.

\dvdelete{Solved: }\pddelete{I think this also should be explained. Why $x_d(p)$ is the probability of the choice? Is this an interpretation? Or this is by construction? DO actually the entries of $x_d(p)$ sum up to 1?}
\dvdelete{Solved: }\pddelete{it is better to explain since the word "stochastic" appears too suddenly. Something like it is too expensive to calculate the whole sum, so we will randomly pick up a component in the sum, which leads to stochastic gradients}

\dv{To use optimization methods for dynamic pricing we can utilize that $f$ is of the form of sum. If the number of suppliers $S$ and consumers $D$ is big, evaluation of all the term in both sums is too expensive in term of working time and algorithmic complexity. But if we randomly pick up only one of $(S+D)$ terms at iteration, computing resources are used much more sparingly, which leads to faster operation of the algorithm. So, considering only one term of sum from \eqref{grad}, we introduce the stochastic gradient oracle:}
\begin{equation}\label{useless_oracle}
\nabla f_i (p) = \begin{cases}
  y_s(p)  & i = s \leq S\\    
  -x_d(p) & i = S + d
\end{cases},\quad i \in \{1, ..., S+D\}
\end{equation}
\dv{Note, that the use of stochastic oracle is very natural. Indeed, using the gradient dynamics of prices from \cite{muller2021dynamic}, we evaluate all of the terms in \eqref{useless_oracle} at iteration, i.e. it is necessary to consider the behaviour of all the consumers and suppliers to make one step of dynamic. But in practice we cannot guarantee that we collect all this information in short time (consumers may be impermanent, so the waiting time may be arbitrarily long). At the same time, the decisions of consumers and suppliers are not rigidly connected: we can observe some number of consumer's sales daily (represented by $x_d(\cdot)$), and much less often and independently the periodic store assortment changes (represented by $y_d(\cdot)$). Therefore, using the dynamics with stochastic oracle in the form of \eqref{useless_oracle} we can immediately take into account newly observed behaviour and make iteration, that is now cheap both in computation and in the required downtime.} However, in real-life environment we alsp cannot estimate the probability of customer's choice, represented by the vector $x_d(p)$ included in second case of \eqref{useless_oracle}. We can obtain only some single sales, those are in fact the random samples of the form $\mathcal{X}_d(p) = (0\;\dotsi 0\;1\;0\;\dotsi 0)$, where $[\mathcal{X}_d(p)]_i = 1$ w.p. $[x_d(p)]_i$, and $[\mathcal{X}_d(p)]_j = 0$ for $i\neq j$. Therefore, $\mathbb{E}[\mathcal{X}_d(p)] = x_d(p)$, and we can introduce another, more practical, stochastic gradient oracle for constructing our dynamics:
\begin{equation}\label{oracle}
\widetilde{\nabla} f_i (p) = \begin{cases}
  y_s(p)  & i = s \leq S\\    
  -\mathcal{X}_d(p) & i = S + d
\end{cases},\quad i \in \{1, ..., S+D\}.
\end{equation}

Now, let us clarify properties of the objective function and introduced oracles. In our simplified setting the following result holds:
\begin{lemma}{(Theorem 3.7 \cite{muller2021dynamic})}\label{lemma}
    $f$ \pd{has $L$-Lipschitz continuous gradient} w.r.t.\dvdelete{Solved: }\pddelete{It is a bit strange that the norm is squared here} $\|\cdot\|_2$ \pd{with}
    \[
        L = \sum_{s=1}^S \frac{1}{\Gamma_s} + \sum_{d=1}^D \frac{1}{\min_j \mu_j},
    \]
    and each $\nabla f_i(p)$ for all $i\in \{1, ..., S+D\}$ is $L_i$-Lipschitz continuous w.r.t. $\|\cdot\|_2$ \pd{with}
    \[
        L_i = 
\begin{cases}
  \displaystyle \frac{1}{\Gamma_i}  & i \leq S\\    
  \displaystyle \frac{1}{\min_j \mu_j} & i > S
\end{cases},\quad i \in \{1, ..., S+D\}
    \]
\end{lemma}

\section{Algorithms and theoretical guarantees}\label{algorithms}

\subsection{Stochastic gradient descent}\label{SGD}

Let us consider the simplest dynamic of prices,
\dvdelete{Solved: }\pddelete{dynamics of what is meant here?}
based on the classical stochastic gradient descent method. A very natural interpretation for this dynamics is that in every iteration we can observe only one of the market participants, supplier or consumer. As soon as participant makes an economical decision (consumer chooses the product or supplier modifies the supply plan), we evaluate the stochastic oracle $\widetilde{\nabla} f_i(p_t)$, and make a step to the equilibrium prices. This dynamic is listed as Algorithm \ref{alg_sgd}. We use the $[\cdot]_+$ notation for the positive part function, i.e. for $a = [b]_+$ we have $a_i = \max\{0, b_i\}$, and denote by $i \sim \mathcal{U}\{1, ..., S+D\}$ the i.i.d. random variables from discrete uniform distribution.

We analyse this dynamic as the \pd{projected} stochastic gradient method with Polyak--Ruppert averaging with tunable parameter $C > 0$ to control the step size of the method. One practical advantage of this method is the robustness to the choice of $C$: it may be chosen regardless of theoretical value of $L$. However, the analysis additionally requires the condition of stochastic gradient's boundedness, but due to the smoothness of $f$ we can bound the norm of \eqref{useless_oracle}, while the additional randomization in the second case of \eqref{oracle} acts on the standard simplex and therefore is also bounded. \dvdelete{I've tried to clarify it. } \pd{Unfortunately, I didn't understand the previous sentence.}
\begin{algorithm}[H]
    \caption{SGD dynamic} \label{alg_sgd}
    \begin{algorithmic}[1]
    \Require $p_0$~--- starting prices values, $N$~--- number of iterations, $C$~--- parameter to control the step size 
    	
    	\For{$t=0,\, 1, \, \ldots, \, N-1$}
    	    \State $i \sim \mathcal{U}\{1, ..., S+D\}$ \dvdelete{Solved: }\pddelete{I suggest to explain the notation for uniform distribution and for the $(\cdot)_+$ function}
    	    
    	    \State $\displaystyle p_{t+1} = \left[p_t - \frac{C}{\sqrt{t+1}} \widetilde{\nabla} f_i (p_t)\right]_+$
        \EndFor
        
        \State $\displaystyle \widetilde{p}_{N} = \frac{1}{N}\sum_{t=1}^N p_t$
    	
    	\State \Return $\widetilde{p}_N$
    \end{algorithmic}
\end{algorithm}
\begin{theorem}{(Theorem 7 \cite{moulines2011non})} \label{conv_sgd}
    Let us assume that used stochastic oracle is uniformly bounded: $\|\widetilde{\nabla} f_i(p_t)\|_2 \leq B$ \pd{for all $i \in \{1, ..., S+D\}$ and $t \in \{1, ..., N\}$}. The suboptimality of prices $\widetilde{p}_N$ given by SGD dynamic (Algorithm \ref{alg_sgd}) is decreasing as follows
    \[
    \mathbb{E}[f(\widetilde{p}_N)] - f(p_*) \leq \frac{\|p_0 - p_*\|_2^2 + C B^2 (1 + C \ln{N})}{2C\sqrt{N}}.
    \]
    Moreover, to obtain the prices satisfying the suboptimality bound
    \[
    \mathbb{E}[f(\widetilde{p}_N)] - f(p_*) \leq \varepsilon,
    \]
    it is \pd{sufficient} to call $\widetilde{\nabla} f_i$ oracle $\displaystyle \mathcal{O}\left(\frac{1}{\varepsilon^2}\right)$ times.
\end{theorem}

Therefore, considered SGD dynamic with Polyak--Ruppert averaging obtains convergence rate of $\displaystyle \mathcal{O}\left(N^{-1/2}\right)$, up to a logarithmic term. This also matches the result from \cite{nesterov2017distributed}.
\dvdelete{Solved above: }\pddelete{I think it is better also to comment on the cheapness of the iteration. If we apply non-stochastic method the rate is much better. So we need to answer the question: what is the benefit of our method.}\dvdelete{Resolved: }\pddelete{Also maybe we can take from the paper \cite{muller2021dynamic} some interpretation of the step. In particular, can the stochastic gradient be interpreted as some price? then it is natural that we subtract a price from a price $p_t$ to adjust the prices.}\dvdelete{In truth, this model is completely irrelevant to dimensional analysis}

\subsection{Adaptive stochastic gradient method}
In this section we describe a slightly different AdaGrad \cite{duchi2011adaptive} dynamics that has a different step size policy. More precisely, the stepsize is chosen based on the stochastic subgradients on the trajectory of the method. This allow the algorithm to adapt to the local information and possibly make longer steps.


\begin{algorithm}[H]
    \caption{AdaGrad dynamic} \label{alg_adagrad}
    \begin{algorithmic}[1]
    \Require $p_0$~--- starting prices values, $N$~--- number of iterations, $\eta$~--- \pd{step size parameter}, $\epsilon$~--- small term to prevent zero division
    	
    	\State $H_0 = 0$
    	
    	\For{$t=0,\, 1, \, \ldots, \, N-1$}
    	    \State $i \sim \mathcal{U}\{1, ..., S+D\}$
    	    \State $g_t = \widetilde{\nabla} f_i (p_t)$
    	    \State $H_{t+1} = H_t + \langle g_t, g_t \rangle$ 
    	    \State $\displaystyle p_{t+1} = \left[p_t - \frac{\eta}{\sqrt{H_{t+1} + \epsilon}} g_t \right]_+$
        \EndFor
        
        \State $\displaystyle \widetilde{p}_{N} = \frac{1}{N}\sum_{t=1}^N p_t$
    	
    	\State \Return $\widetilde{p}_N$
    \end{algorithmic}
\end{algorithm}
\pddelete{I think, the following theoretical result is proved in \cite{duchi2018introductory} for a different algorithm} \dvdelete{I'm not sure, the $\sqrt{H}/\eta$ is exactly the $H$ in (4.3.7), it seems to be classical simplified adagrad}
\begin{theorem}{(Corollary 4.3.8 \cite{duchi2018introductory})} \label{conv_adagrad}
    Let $\|p_t - p_0\|_\infty \leq R$ for all $t \in \{1, ..., N\}$, $\eta$ is proportional to $R$\dvdelete{  Solved: }\pddelete{it is better to explain what is meant by $\sim$ here}. The suboptimality of the prices $\widetilde{p}_N$ given by AdaGrad dynamic (Algorithm \ref{alg_adagrad}) is decreasing as follows
    \[
    \mathbb{E}[f(\widetilde{p}_N)] - f(p_*) \leq \frac{3 R}{2 N} \cdot \sum_{i=1}^n \mathbb{E} \left[\left(\sum_{t=1}^N [g_t]_i^2\right)^{1/2} \right].
    \]
    To obtain the prices satisfying the suboptimality bound
    \[
    \mathbb{E}[f(\widetilde{p}_N)] - f(p_*) \leq \varepsilon,
    \]
    it is \pd{sufficient} to call $\widetilde{\nabla} f_i$ oracle $\displaystyle \mathcal{O}\left(\frac{1}{\varepsilon^\dv{2}}\right)$\dvdelete{ Of course. Solved: }\pddelete{shouldn't here $\varepsilon$ be squared?} times.
\end{theorem}

Hence the proposed AdaGrad dynamic demonstrates convergence rate of the order $\displaystyle \mathcal{O}\left(N^{-1/2}\right)$\dvdelete{Of course. Solved: }\pddelete{I'm not sure about this. I think, the sum of $g_t$ grows like $\sqrt{N}$}. 
This asymptotic is similar to that for SGD dynamic, but in practice such a simple modification of the step size allows to discernibly improve the convergence rate.
At the same time, the step size hyperparameter $\eta$ is still free and it allows to manually tune the algorithm for the best practical efficiency.

\section{The case of infinite number of consumers}\label{infty_d}

In general, the total expected revenue framework considered in \cite{muller2021dynamic} allows one to describe not only the individual consumers, but also the groups of consumers with \pd{close} behavior. At the same time, their model covers only the setting of finite number of agents, which may be not completely practical. Indeed, if the number of agents is huge or new agents can enter the marketplace as time goes, it makes sense to consider the limit when the number of agents tends to infinity. This leads to a general, non-finite-sum objective given as an expectation. 

To be more specific, we consider that consumers are represented by a random variable with some unknown distribution, i.e. $d \sim \mathcal{D}$. Then the characteristic vector $(a_{i d})_{i=1}^n$ is also random. To maintain the property of market clearing at optimal prices we also assume that the number of suppliers is infinite and that they are represented by another random variable, i.e. $s \sim \mathcal{S}$. The next step is to take the limit in \eqref{finite_problem} and make transition to expectation w.r.t. distributions $\mathcal{D},\mathcal{S}$ instead of sums. For the sake of normalization we introduce the parameter $0 < \beta < 1$ that is equal to the fraction of suppliers among all market participants. Informally, $\beta = \lim_{(S+D) \rightarrow \infty} S / (S + D)$.

In this way, considering function $\frac{1}{S+D} f(p)$ instead of $f(p)$ and taking the limit as $(S+D) \rightarrow \infty$, we have the new optimization problem in the form of expectation:
\begin{equation}\label{infinite_problem}
    \min_{p \in \mathbb{R}^n_+} \left\{ \widetilde{f}(p) := \beta \cdot \mathbb{E}_{s \sim \mathcal{S}} [\pi_s(p)] + (1 - \beta) \cdot \mathbb{E}_{d \sim \mathcal{D}} [E_d(p)] \right\}.
\end{equation}

Note, that to preserve the convergence properties of the SGD dynamic considered in Section \ref{SGD} it is sufficient just to generalize the used stochastic oracle \eqref{oracle} to the proposed setting by defining
\begin{equation}\label{infinite_oracle}
\widetilde{\nabla} \widetilde{f}(p) := \begin{cases}
  y_s(p)  & \text{for }s \sim \mathcal{S}\text{ w.p. }\beta \\    
  -\mathcal{X}_d(p) & \text{for }d \sim \mathcal{D}\text{ w.p. }1-\beta
\end{cases}.
\end{equation}
\pd{Thus, to generate the stochastic gradient, we first with probability $\beta$ choose to choose among suppliers or with probability $1-\beta$ we choose to choose among consumers. Then, in the former case we sample supplier from the distribution  $\mathcal{S}$ and in the latter case we sample consumer from the distribution  $\mathcal{D}$. Finally, for the chosen agent the stochastic gradient is defined as in \eqref{oracle}.}

Since problem \eqref{infinite_problem} is a general stochastic optimization problem, we can apply Algorithm \ref{alg_sgd} and obtain the dynamic listed as Algorithm \ref{alg_sgd_online}. \dv{Its interpretation is quite similar to that of Algorithm \ref{alg_sgd}: at the every iteration we observe the behaviour of one market participant, consumer or supplier, and change the prices in a proper way. As in the SGD dynamic, samplings $d \sim \mathcal{D}$, $s \sim \mathcal{S}$ and switching w.p. $\beta$ are provided by natural flow of participants, we assume that information about participants decisions arrives uniformly.} \dvdelete{Solved:} \pddelete{Unfortunately I didn't understand the previous sentence.} Since Lemma \ref{lemma} and the conditions of Theorem \ref{conv_sgd} still hold, the convergence rate of dynamic below is similar to the one given in Theorem \ref{conv_sgd}.
\begin{algorithm}[H]
    \caption{SGD dynamic (online setting)} \label{alg_sgd_online}
    \begin{algorithmic}[1]
    \Require $\beta$~--- fraction of suppliers, $p_0$~--- starting prices values, $N$~--- number of iterations, $C$~--- parameter to control learning rate
    	
    	\For{$t=0,\, 1, \, \ldots, \, N-1$}
    	    \State $\displaystyle p_{t+1} = \left[p_t - \frac{C}{\sqrt{t+1}} \widetilde{\nabla} \widetilde{f} (p_t)\right]_+$
        \EndFor
        
        \State $\displaystyle \widetilde{p}_{N} = \frac{1}{N}\sum_{t=1}^N p_t$
    	
    	\State \Return $\widetilde{p}_N$
    \end{algorithmic}
\end{algorithm}

\section{The case of zero quantity adjustment costs}\label{zero_gamma}
\dvdelete{Solved:} \pddelete{I think, it is better to find a title without formulas}

In this section we return to the setting of Section \ref{S:problem} and consider the limiting case when in  \eqref{pi} $\Gamma_s=0$. In this case there is no typical supply $\hat{y}_s$ and all the necessary information about the costs incurred by the supplier is set by the function $c_s(\cdot)$. In this case there is no guarantee that the function $f$ in \eqref{finite_problem} is smooth (cf. Lemma \ref{lemma}) and the algorithms have to be properly modified to guarantee the convergence. This section is devoted to such modifications.

\subsection{Mirror descent dynamic}\label{SMD}

Stochastic mirror descent (SMD) \cite{duchi2012ergodic,nemirovski2009robust,nemirovskij1983problem} is widely used and theoretically optimal algorithm for stochastic convex non-smooth optimization problems. We consider here a particular case of the SMD dynamic with Euclidean proximal setup, in which it looks quite similar to SGD dynamic, but uses a different step size policy. SMD dynamic is listed below as Algorithm \ref{alg_smd}.
\begin{algorithm}[H]
    \caption{SMD dynamic} \label{alg_smd}
    \begin{algorithmic}[1]
    \Require $p_0$~--- starting prices values, $N$~--- number of iterations, $C$~--- parameter to control step size
    	
    	\For{$t=0,\, 1, \, \ldots, \, N-1$}
    	    \State $i \sim \mathcal{U}\{1, ..., S+D\}$
    	    
    	    \State $\displaystyle p_{t+1} = \left[p_t - \frac{C R}{M \sqrt{t+1}} \widetilde{\nabla} f_i (p_t)\right]_+$
        \EndFor
        
        \State $\displaystyle \widetilde{p}_{N} = \frac{1}{N}\sum_{t=1}^N p_t$
    	
    	\State \Return $\widetilde{p}_N$
    \end{algorithmic}
\end{algorithm}

Theoretical guarantee below utilizes the Lipschitz continuity of function $f$, more precisely the boundedness of the stochastic gradient norm. Considering the expression \eqref{grad}, we can bound the first term by applying some economic reasoning, namely that the supply is always limited or scarcity principle \cite{burke1800thoughts}. The second term is also bounded since all $x_d$ are probability vectors and belong to the standard simplex. Stochastic gradient given by \eqref{oracle} is bounded by the same reason. 

\begin{theorem}{(Proposition 1 \cite{lan2012validation})} \label{conv_smd}
    Let $\|p_t - p_0\|_2 \leq R$ and $\mathbb{E}[\|\widetilde{\nabla}f(p_t)\|^2] \leq M^2$ for all $t \in \{1, ..., N\}$. The suboptimality of prices $\widetilde{p}_N$ given by SMD dynamic (Algorithm \ref{alg_smd}) is decreasing as follows
    \[
    \mathbb{E}[f(\widetilde{p}_N)] - f(p_*) \leq \frac{\max\{C, C^{-1}\} R M}{\sqrt{N}}.
    \]
    Moreover, to obtain the prices satisfying the suboptimality bound
    \[
    \mathbb{E}[f(\widetilde{p}_N)] - f(p_*) \leq \varepsilon,
    \]
    it is \pd{sufficient} to call $\widetilde{\nabla} f_i$ oracle $\displaystyle \mathcal{O}\left(\frac{1}{\varepsilon^2}\right)$ times.
\end{theorem}

So, the described SMD dynamic has convergence rate of $\displaystyle \mathcal{O}\left(N^{-1/2}\right)$ in non-smooth case that takes place if $\Gamma_s = 0$ for some $s$. This bound matches that of the SGD dynamic from Theorem \ref{conv_sgd}. 

\dvdelete{Solved:}
\pddelete{Unfortunately, the previous two sentences I didn't understand.}

\subsection{Dual smoothing}
The special structure of the function $\pi_s(p)$ in \eqref{pi} allows us to use the Nesterov's smoothing technique \cite{nesterov2005smooth}. The idea is, in the case when $\Gamma_s = 0$, to replace the $- \Gamma_s \cdot \|y - \hat{y}_s\|_2^2$ term in \eqref{pi} with the synthetic penalty $-\eta \cdot \|y - y_0\|_2^2$, where $\eta = \varepsilon / (2 R^2)$ for some fixed target suboptimality  $\varepsilon$ and $R$ such that $\|y_* - y_0\|_2 \leq R$. With this substitution for all $s$, following the argumentation of Lemma \ref{lemma}, we have that the modified function $f_{\eta}$ has Lipschitz-continuous gradient with constant
\[
L = \frac{2 S R^2}{\varepsilon} + \sum_{d=1}^D \frac{1}{\min_j \mu_j}.
\]
At the same time, if we minimize $f_{\eta}$ up to accuracy $\frac{\varepsilon}{2}$, i.e. we have for some $\widetilde{p}$ that
\[
f_{\eta}(\widetilde{p}) - \min_{p \in \mathbb{R}^n_+}f_{\eta}(p) \leq \frac{\varepsilon}{2},
\]
then it holds that
\[
f(\widetilde{p}) - f(p_*) \leq \varepsilon.
\]
It means that we can obtain the solution of the problem \eqref{finite_problem} satisfying the target suboptimality bound by optimizing the modified function $f_{\eta}$ that is smooth and therefore allows us to apply some of the methods described in previous sections. 

Note that the transition from non-smooth setting to the smooth one with the described approach is not free in terms of convergence rate due to the dependence of the constant $L$ on the target suboptimality $\varepsilon$. 

\section{Numerical experiments}\label{experiments}

In this section, we focus on the problem \eqref{finite_problem}, which is motivated by important applications to offering smart pricing options by online marketplaces and management of demand and supply by such financial intermediaries like brokers \cite{muller2021dynamic}. 

We provide the results of our numerical experiments, which are performed on a PC with processor Intel Core i7-8650U 1.9 GHz using pure Python 3.7.3 (without C code) under managing OS Windows 10 (64-bits). Numpy.float128 data type with precision $1e-18$ and with max element $\approx 1.19e+4932$ is used. Random seed is set to 17. 

We compare Algorithm \ref{alg_sgd} and Algorithm \ref{alg_adagrad} with pricing dynamics 4.2 and 4.4 from \cite{muller2021dynamic} on the problem \eqref{finite_problem} with the following settings: number of suppliers $S = 5$, number of consumers $D = 10$, number of products $n = 20$, number of groups $m = 5$. $\Gamma_s$ is chosen like $10^{-4}$ (some small value that affects to convergence of pricing schemes from \cite{muller2021dynamic}). We generated also $S$ vectors $\hat{y}_s$ from uniform distribution $\mathcal{U}[0.01, 2]$ of size $n$. Initialization point $p$ is chosen from uniform distribution $\mathcal{U}[0.01, 5]$ of size $n$, but it is scaled by the maximum element. Matrix $A$ consists of columns with each one is drawn from uniform distribution
And $\mu$ values are drawn i.i.d. from uniform distribution $\mathcal{U}[0.1, 1]$ of size $m$.

We consider $c_s(y)$ equal to $||y||_2^2$. Hence, the closed form solution of \eqref{pi} is
\begin{equation*}
y_s(p) = \frac{p + 2 \Gamma_s \hat{y}_s}{2 (1 + \Gamma_s)}.
\end{equation*}
To estimate the suboptimality we perform pricing dynamic 4.2 \cite{muller2021dynamic} (Gradient descent method) with stopping criterion $||p_{t+1} - p_{t}||_2 \leq 10^{-10}$ and obtain lower bound for $f_* = f(p_*)$ from Theorem 4.3 \cite{muller2021dynamic}. 

\begin{figure}[H]
	\centering
	\vspace{-0.3cm}
	\includegraphics[width=0.70\textwidth]{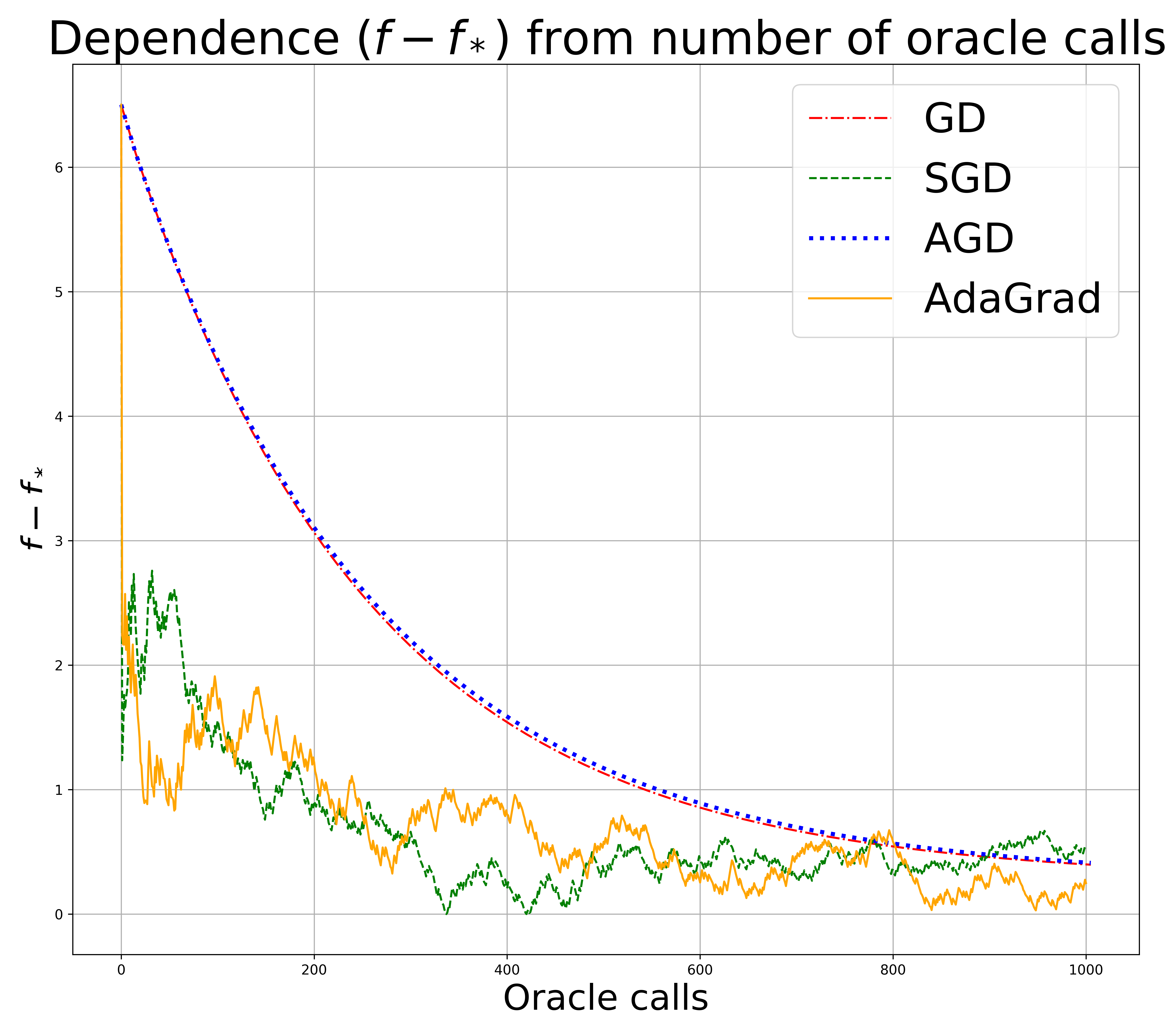}
	\vspace{-0.4cm}
	\caption{Dependence of $f-f_*$ from the number of oracle calls.}
	\label{fig:dependence_from_oracle_calls}
\end{figure}

During the experiments we store objective suboptimality and number of oracle calls for each algorithm at every iteration. The results are presented on the Figure~\ref{fig:dependence_from_oracle_calls}.  

As we can see from the Figure \ref{fig:dependence_from_oracle_calls}, there is almost no difference between accelerated and non-accelerated pricing schemes from \cite{muller2021dynamic} and that the stochastic algorithms converge to the optimal value much faster than accelerated gradient method with respect to chosen parameters.

\section*{Conclusion}

We propose a stochastic version of the formulation of the problem of finding equilibrium prices by minimizing the potential function of total expected revenue, proposed in a deterministic form in \cite{muller2021dynamic}. Thanks to the introduction of the stochastic oracle, we analyze and interpret in terms of observing the real-life marketplace environment several dynamics based on the efficient stochastic gradient optimization methods. We also analyze the case of a non-smooth potential function in the case of zero quantity adjustment costs. In addition, we propose a generalized setting allowing an infinite number of market participants and equally suitable application of optimization methods. Numerical experiments show that stochastic methods turn out to be discernibly more efficient in comparison with full-gradient dynamics.

\bibliographystyle{splncs04}
\bibliography{mybibliography}
\end{document}